\documentclass{article}
\usepackage{amssymb}

\begin{document}
\begin{center}
{\large\bf Pointwise convergence of solution to Schr\"{o}dinger equation on manifolds}
\end{center}
\begin{center}
Xing Wang\footnote{Department of Mathematics, Wayne State University, Detroit, Michigan, 48202, USA}, Chunjie Zhang\footnote{
Corresponding author, Department of Mathematics, Hangzhou Dianzi University, Hangzhou, 310018, China }
\end{center}
{\bf Abstract.} Let $(M^n,g)$ be a Riemannian manifold without boundary. We study the amount of initial regularity is required so that the solution to free Schr\"{o}dinger equation converges pointwisely to its initial data. Assume the initial data is in $H^\alpha(M)$. For Hyperbolic Space, standard Sphere and the 2 dimensional Torus, we prove that $\alpha>\frac{1}{2}$ is enough. For general compact manifolds, due to lacking of local smoothing effect, it is hard to beat the bound $\alpha>1$ from interpolation. We managed to go below 1 for dimension $\leq 3$. The more interesting thing is that, for 1 dimensional compact manifold, $\alpha>\frac{1}{3}$ is sufficient.


\section{Introduction}
In the Euclidean setting, Carleson \cite{Ca} proposed a question regarding the amount of regularity required on the initial data $f$ so that
$$e^{-it\triangle}f(x)=\int_{\Bbb R^n}e^{i(\langle x,\xi\rangle+t|\xi|^2)}\hat{f}(\xi)d\xi\rightarrow f(x),~~{\mbox a.e.}~x\in \Bbb R^n$$
as $t$ goes to $0$. Here $e^{-it\triangle }f(x)$ is the solution to free Schr\"{o}dinger equation
$$\left\{\begin{array}{ll}
i\partial_tu-\Delta u=0,& (t,x)\in \Bbb R\times \Bbb R^n\\
u(0,x)=f(x),& x\in \Bbb R^n.\end{array}\right.$$
The problem then has been treated by many authors. When $n=1$, Carleson himself proved that $f\in H^{\frac{1}{4}}(\Bbb R)$ is sufficient. B.E.J. Dahlberg and C.E. Kenig \cite{DK} proved this is necessary for all dimensions.  For higher dimensions, M. Cowling \cite{Cow} studied this problem for a general class of self-adjoint operators and obtained $\alpha>1$ for the Schr\"{o}dinger operator. Later on, P. Sj\"{o}lin \cite{Sj} proved a local smoothing effect and thus improved the bound to $\alpha>\frac{1}{2}$, which was also independently proved by L. Vega in \cite{Ve1}.

For n=2, J. Bourgain \cite{B1} showed we can go below $\frac{1}{2}$ a little bit. Furthermore, continuous improvement has been made by Moyua, Vargas, Vega \cite{MVV}, Tao, Vargas
\cite{TV1,TV2}. Right now, the best result known for n = 2 is $f\in H^\alpha$, $\alpha>\frac{3}{8}$ by S. Lee (\cite{Lee}).

For $n\geq 3$, Bourgain\cite{B2} beat the bound $\alpha>\frac{1}{2}$ and showed that the solution to free Schr\"{o}dinger equation converges pointwisely to its initial data $f$, provided $f\in H^\alpha(\Bbb R^n)$ where $\alpha>\frac{1}{2}-\frac{1}{4n}$. Furthermore, he also showed that
when $n>4$, $\alpha>\frac{1}{2}-\frac{1}{n}$ is necessary. In \cite{LR}, Luc$\acute{a}$ and Rogers refined the necessary condition to $\alpha>\frac{1}{2}-\frac{1}{n+2}$.

In this paper, we deal with a similar problem in the manifold setting. Below we always take $(M^n,g)$ to be a complete manifold endowed with
a $C^\infty$ metric $g$. Denote by $\Delta $ the Laplace-Beltrami operator associated with $g$.

The  Schr\"{o}dinger equation on $(M^n,g)$ is given by
$$\left\{\begin{array}{ll}
i\partial_tu-\Delta u=0,& (t,x)\in \Bbb R\times M,\\
u(0,x)=f(x),& x\in M.\end{array}\right.\eqno(1.1)$$

Specifically, if $M$ is compact, it has discrete spectrum, and there exists an orthonormal basis
$\{ e_j\}$ of eigenfunctions such that
$$\Delta e_j(x)=\lambda_j^2 e_j, \quad e_j\in C^\infty(M), \quad
\int_{\mathcal{M}}e_je_k\, dV_g=\delta_{jk}.    $$

The
eigenvalues $ 0=\lambda_0<\lambda_1\leq \lambda_2\leq
\lambda_3,\cdots, $ are listed in ascending order and counted by
multiplicity.

We can then write the solution of (1.1) as
$$e^{-it\triangle }f(x)=\sum\limits_{j=0}^\infty e^{it\lambda_j^2}\hat{f}_je_j(x),\eqno(1.2)$$
where $\hat{f}_j$ is the j-th Fourier coefficient given by
$$\hat{f}_j=\int_M f(x)e_j(x)dV_g.$$

When $M=\Bbb H^n, \Bbb S^n$ or $\Bbb T^2$, we are able to obtain the same result as in P. Sj\"{o}lin's work \cite{Sj} for $\Bbb R^n$.\\\\
{\bf Theorem 1.1.} The solution $u(t,x)$ to equation (1.1) converges pointwisely to the initial data $f$, whenever $f\in H^\alpha(M), \alpha>\frac{1}{2}$. Here  $M=\Bbb H^n, \Bbb S^n$ or $\Bbb T^2$ endowed with standard metric.\\\\
{\bf Remark 1.1.} These three cases are actually proved via three different methods. For Hyperbolic Space, we derive it by showing a local smoothing effect. For standard Sphere, we take full advantage of the spectrum concentration. For the 2 dimensional Torus, we apply our argument for the general manifold and combine it with Strichartze estimates on 2 dimensional Torus obtained by Bourgain \cite{B4}.

For general compact manifold, if $\alpha>1$, there is a way to prove the pointwise convergence quickly, using a method from \cite{BG}. In fact by the
expression (1.2) and Parseval's formula, we easily have
$$\|e^{-it\Delta}f(x)\|_{L^2(M,L^2(0,1])}\leq \|f\|_{L^2(M)},$$
$$\|e^{-it\Delta}f(x)\|_{L^2(M,H^1(0,1])}=\|(e^{-it\Delta}-i\Delta e^{-it\Delta})f(x)\|_{L^2(M,L^2(0,1])} \leq \|f\|_{H^2(M)}.$$
Interpolating between the two yields
$$\|e^{-it\Delta}f(x)\|_{L^2(M,H^s(0,1])}\leq \|f\|_{H^{2s}(M)},$$
which, when combined with the Sobolev imbedding, further leads to
$$\left\|\sup_{0<t\leq 1}|e^{-it\Delta}f(x)|\right\| \lesssim \|f\|_{H^{2s}(M)},~s>\frac{1}{2}.\eqno(1.3)$$
The pointwise convergence then follows from a standard argument using the above boundedness.

From here and below, we will use the notation $A\lesssim B$
to mean that there is some constant $C$ independent of all essential variables, such that $A\leq C B$.

For convenience, we denote
$$T^*f(x)=\sup_{0<t\leq 1}|e^{-it\Delta}f(x)|.$$

Inequality (1.3) says that the Maximal Schr\"{o}dinger operator $T^*$ is bounded from $H^\alpha(M),~\alpha>1$ to $L^2(M)$. Then using a standard argument, pointwise convergence follows. From Theorem 1, it is reasonable to conjecture that $\alpha>\frac{1}{2}$ should be sufficient, but it is hard to break the bound $\alpha>1$. This is due to the absence of a local smoothing effect \cite{SD1,SD2}, and we no longer have scaling invariances as in Euclidean Spaces. Fortunately, by utilizing the Strichartz estimate \cite{BGT}, we manage to overcome these difficulties and break the bound $\alpha>1$ in lower dimensions.\\\\
{\bf Theorem 1.2.} Let $(M,g)$ be a connected, compact manifold without boundary of dimension $n$. The solution $u(t,x)$ to equation (1.1) converges pointwisely to the initial data f, whenever $f\in  H^\alpha(M)$, $\alpha>\frac{1}{3}$ for $n=1$, $\alpha > \frac{3}{4}$ for $n = 2$, or $ \alpha> \frac{9}{10}$ for $n = 3$.\\

We now give a brief outline of what follows.
In Section 2, we provide some basic facts about Hyperbolic Space and Sphere. In Section 3, we prove the Hyperbolic Space case in Theorem 1.1. In Section 4, we prove the standard Sphere case, which exhibits some difference between the Hyperbolic Space and the Euclidean case. In Section 5, we prove Theorem 1.2 for n=2,3. In Section 6, we apply the previous results to the Torus case, and several other examples.

\section{Preliminaries}
In this section, we provide some basic facts we will need in the later context.\\
{\bf Definition.} We define the Hyperbolic space as given by the polar parametrization:
$$\Bbb H^n=\{(t,x)\in\Bbb R^{n+1}, (t,x)=(\cosh r,\sinh r\omega),~r>0,~\omega\in \Bbb S^{n-1} \}.$$\\
The metric tensor is given by
$$g=dr^2+\sinh rd\omega^2,$$
where $d\omega^2$ is the standard metric on Sphere.
The volume element is
$$dV=\sinh r^{n-1}drd\omega.$$
The Laplace-Beltrami operator is
$$\Delta_{\Bbb H^n}=\partial_r^2+(n-1)\frac{\cosh r}{\sinh r}\partial_r+\frac{1}{\sinh^2 r}\Delta_{\Bbb S^{n-1}}.$$
The Sobolev space is defined by
$$H^s(\Bbb H^n)=\{f|(1-\Delta_{\Bbb H^n})^{s/2}f\in L^2(\Bbb H^n)\}.$$
When $s=m$ is a positive integer, this is equivalent to the usual definition.
Let $f$ be the function on $\Bbb H^n$ given by
$$F=\sqrt{1+r^2}\eqno(2.1)$$
written in polar coordinate. Noticing that $F\in C^\infty(\Bbb H^n)$, with some direct calculation, one has\\\\
{\bf Proposition 2.1.} Let $F$ be as above. Then
\begin{enumerate}
\item $Hess_F=\frac{1}{\sqrt{1+r^2}}dr\otimes dr+\frac{r}{\sqrt{1+r^2}}Hess_r\geq A(r)g$, where $A(r)=\min(\frac{1}{\sqrt{1+r^2}},\frac{r}{\sqrt{1+r^2}}\frac{\cosh r}{\sinh r})$.
\item $\|\nabla^k F\|\in L^\infty(\Bbb H^n)$ for any integer $k>0$.\\
\end{enumerate}

Now we describe the spectrum concentration of Sphere.\\\\
{\bf Proposition 2.2.}  Denote $\mu_k$ be the $k$-th eigenvalue on standard sphere, without counting multiplicity. Let $E_k$ be the corresponding eigenspace. Then
$\mu_k=k(k+n-1)$ and $\dim E_k\approx k^{n-1}$.\\

In the general compact manifold case, we will use H\"{o}rmander's oscillatory integral estimates, see \cite{SoFIO} Theorem 2.2.1 for a detailed reference.
Consider oscillatory integrals of the form
$$T_hf(z)=\int e^{\frac{i}{h}\phi(z,y)}a(z,y)f(y)dy.$$
Here $a\in C_0^\infty(\Bbb R^{n+1},\Bbb R^n)$, $\phi$ is real and $C^\infty$ in a neighborhood of $supp~a$.
Then the canonical relation associated to $\phi$ is defined as
$$C_\phi=\{(z,\phi_z(z,y),y-\phi_y(z,y))\}\subset T^*\Bbb R^{n+1}\times T^*\Bbb R^{n}.$$
\\
{\bf Lemma 2.1.}  Denote $\Pi_{T^*\Bbb R^{n}}:C_\phi\rightarrow T^*\Bbb R^{n}$ as the natural projection, similarly as $\Pi_{T^*\Bbb R^{n+1}}$. Assume
\begin{enumerate}
\item rank$d\Pi_{T^*\Bbb R^{n}}\equiv 2n;$
\item $S_{z_0}=\Pi_{T_{z_0}^*\Bbb R^{n+1}}C_\phi$ has everywhere non-vanishing Gaussian curvature for any $z_0\in supp_z~a$.
\end{enumerate}
Then
$$\|T_hf\|_{L^p(\Bbb R^{n+1})}\lesssim h^{(n+1)/q}\|f\|_{L^p(\Bbb R^n)}$$
if $q=\frac{n+2}{n}p'$ and $1\leq p\leq 2$ for $n\geq 2$, or $1\leq p<4$ for $n= 1$.

\section{Solution on Hyperbolic Space}
We start with the following local smoothing lemma.\\\\
{\bf Lemma 3.1.} Let $u(t)=u(x,t)$ be the solution of (1.1) with $M=\Bbb H^n$. Denote $B_R$ as the geodesic ball centered at the origin. Then there exists a constant $C=C(n,R)>0$, such that
$$\|u\|_{L^2([0,1]\times B_R)}\leq C\|f\|_{H^{-\frac{1}{2}}(\Bbb H^n)}.\eqno(3.1)$$\\
{\bf Remark 3.1.} The smoothing effect of Schr\"{o}dinger evolution group has been intensively studied. Here we refer the readers to the works \cite{SD1,SD2}. The proof of Lemma 3.1 below follows from \cite{SD1} with some modifications.\\\\
{\bf Proof of Lemma 3.1.} In this section, we will denote $\Delta_{\Bbb H^n}=\Delta$.

Choose $\phi\in C^\infty(\Bbb H^n)$ to be a cut off function such that $\phi\equiv 1$ in $B_R$ and $\phi\equiv 0$ outside $B_{2R}$, $0\leq\phi\leq 1$ and $|\nabla\phi|\lesssim \frac{1}{R}$. Let $f$ be the function defined by (2.1).
Consider the self-adjoint linear differential operator
$$X=\frac{\nabla F}{i}+(\frac{\nabla F}{i})*=\frac{2\nabla F}{i}+\frac{\Delta F}{i}.$$
In local coordinates,
$$\nabla F=g^{ij}\frac{\partial F}{\partial x^i}\frac{\partial}{\partial x^j}.$$
Let $N=1-\Delta$, then $P=N^{-1/4}XN^{-1/4}$ is a pseudodifferential operator of order 0. Since $\|\nabla^k F\|\in L^\infty(\Bbb H^n)$, one can see $X:H^s(\Bbb H^n)\rightarrow H^{s-1}(\Bbb H^n)$ is continuous. Thus $P$ is bounded on $L^2(\Bbb H^n)$.

\begin{eqnarray*}
\frac{d}{dt}(-Pu(t),u(t))&=&(-iP\Delta u(t),u(t))+(-Pu(t),i\Delta u(t))\\
&=&2\Re(-iP\Delta u(t),u(t))\\
&=&2\Re(-i\Delta N^{-1/4}u(t),XN^{-1/4}u(t))\\
&=&\int_{\Bbb H^n}(4Hess_F(\nabla v(t),\overline{\nabla v(t)})-\Delta^2 F |v(t)|^2)dV\\
&\geq & A_{2R}\int_{B_{2R}}\phi|\nabla v(t)|^2- \int_{\Bbb H^n}\Delta^2 F |v(t)|^2dV.
\end{eqnarray*}
Here $v(t)=N^{-1/4}u(t)$ and $A_{2R}=4\inf_{0\leq r\leq 2R}A(r)$. Integrate the above inequality from 0 to 1, we have
\begin{eqnarray*}
\int_{[0,1]\times B_{2R}}\phi|\nabla v(t)|^2dVdt&\lesssim& \|u(0)\|_{L^2(\Bbb H^n)}+\|u(1)\|_{L^2(\Bbb H^n)}+\|v(t)\|_{L^2([0,1]\times\Bbb H^n)}\\
&\lesssim&\|f\|_{L^2(\Bbb H^n)}
\end{eqnarray*}
with a constant depending on $R$.

Notice $v(t)=N^{-1/4}u(t)$ solves the following Schr\"{o}dinger equation:
$$\left\{\begin{array}{ll}
i\partial_tv-\Delta v=0,& (t,x)\in \Bbb R\times \Bbb H^n,\\
v(0,x)=N^{-1/4}f(x),& x\in \Bbb H^n.\end{array}\right.$$

Applying the above inequality one obtains
$$\int_{[0,1]\times B_{2R}}\phi|\nabla N^{-1/2}u(t)|^2dVdt\lesssim \|f\|_{H^{-1/2}(\Bbb H^n)}.$$
Since
\begin{eqnarray*}
(\phi\nabla N^{-1/2}u,\nabla N^{-1/2}u)&=&(\phi(-\Delta) N^{-1/2}u,N^{-1/2}u)+(\nabla\cdot\phi\nabla(N^{-1/2}u),N^{-1/2}u)\\
&=&(N^{1/2}u, \phi N^{-1/2}u)-(\phi N^{-1/2}u,N^{-1//2}u)\\
&&+(\nabla\phi\cdot\nabla(N^{-1/2}u),N^{-1/2}u)\\
&\geq&(u, N^{1/2}(\phi N^{-1/2}u))-C_1\|u\|_{H^{-1/2}(\Bbb H^n)}
\end{eqnarray*}

By the sharp G$\mathring{a}$rding Inequality(see \cite{Ho}), we have
$$(u, N^{1/2}(\phi N^{-1/2}u))\geq (\phi u,u)-C_2\|u\|_{H^{-1/2}(\Bbb H^n)}.$$
Here $C_1$, $C_2$ are constants depending on $R$.

Combining the estimates above, we finish the proof of the lemma.
As a corollary, we have\\\\
{\bf Corollary 3.1.} The following estimate holds
$$\|\Delta u\|_{L^2([0,1]\times B_R)}\leq C(R)\|f\|_{H^{3/2}(\Bbb H^n)}.\eqno(3.2)$$\\
{\bf Proof of Theorem 1.1 for $\Bbb H^n$.} The proof is similar to the argument for (1.3) in the introduction, we finish the proof by interpolation between (3.1) and (3.2).\\\\
{\bf Remark 3.1.} The argument above also applies to other manifolds which the local smoothing effect holds. As a consequence, if we perturb the the standard metric of Euclidean or Hyperbolic space in a finite domain, such that there's no trapped geodesics, then the local smoothing effect still holds. See \cite{SD1,SD2} for more examples.

\section{Solution on Sphere}
{\bf Proof.} By Proposition 2.2, we know that the $k$-th eigenvalue of $-\Delta_{\Bbb S^n}$ is $\mu_k=k(k+n-1)$, and
the the eigenfunctions attached to $\mu_k$ are the sphere harmonics of degree $k$, which form a linear space of dimension $d_k\approx k^{n-1}$.
Take $e_{j_1}(x),~e_{j_2}(x),\cdots,e_{j_{d_k}}(x)$ to be an $L^2$ normalized base of this linear space. Then for each $f\in C^\infty(\Bbb S^n)$, we have
$$e^{-it\Delta_{\Bbb S^n}}f(x)=\sum_{k=0}^{+\infty}\sum_{l=1}^{d_k}e^{-it\mu_k}\hat{f}_{k_l}e_{k_l}(x).\eqno(4.1)$$
We wish to prove
$$\|\sup_{0<t\leq 1}|e^{-it\Delta_{\Bbb S^n}}f(x)|\|_{L^2(\Bbb S^n))}\lesssim \|f\|_{H^\alpha(\Bbb S^n)},~\alpha>\frac{1}{2}.$$
It suffices to bound $e^{-it(x)\Delta_{\Bbb S^n}}f$, and for this, we have by (4.1),
\begin{eqnarray*}
\|e^{-it(x)\Delta_{\Bbb S^n}}f(x)\|_{L^2}&=& \left\|\sum_{k=0}^{+\infty}e^{-it(x)\mu_k}\sum_{l=1}^{d_k}\hat{f}_{k_l}e_{k_l}(x)\right\|_{L^2}\\
&\lesssim & \sum_{k=0}^{+\infty}\left\|\sum_{l=1}^{d_k}\hat{f}_{k_l}e_{k_l}(x)\right\|_{L^2}\\
&\lesssim & \sum_{k=0}^{+\infty}\left(\sum_{l=1}^{d_k}|\hat{f}_{k_l}|^2\right)^{1/2}\\
&\lesssim & \sum_{k=0}^{+\infty}(1+\mu_k)^{-\alpha/2}\left(\sum_{l=1}^{d_k}(1+\mu_k)^{\alpha}|\hat{f}_{k_l}|^2\right)^{1/2}\\
&\lesssim & \left(\sum_{k=0}^{+\infty}(1+\mu_k)^{-\alpha}\right)^{1/2}\left(\sum_{k=0}^{+\infty}\sum_{l=1}^{d_k}(1+\mu_k)^{\alpha}|\hat{f}_{k_l}|^2\right)^{1/2}\\
&=&\left(\sum_{k=0}^{+\infty}(1+k(k+n-1))^{-\alpha}\right)^{1/2}\|f\|_{H^\alpha(\Bbb S^n)}\\
&\leq & C_\alpha \|f\|_{H^\alpha(\Bbb S^n)}.
\end{eqnarray*}

\section{General case}
To prove Theorem 1.2, we first do a spectrum decomposition. Take $\tilde{\psi}\in C^\infty_0(\Bbb R)$ and $\psi\in C_0^\infty(\Bbb R\backslash\{0\})$ such that
$$\tilde{\psi}(\lambda^2)+\sum\limits_{k=1}^{+\infty}\psi(2^{-2k}\lambda^2)=1,~~\forall~\lambda.$$
Then we have the decomposition
$$f=\tilde{\psi}(\Delta^2)f+\sum\limits_{k=1}^{+\infty}\psi(2^{-2k}\Delta^2)f$$
for any $f\in C^\infty(M)$, and furthermore,
$$T^*f\lesssim T^*(\tilde{\psi}(\Delta^2)f)+\sum\limits_{k=1}^{+\infty}T^*(\psi(2^{-2k}\Delta^2)f).\eqno(5.1)$$
As mentioned in the introduction, we only need to show
$$\|T^*f\|_{L^p(M)}\lesssim \|f\|_{H^\alpha(M)},~~p=2(n+2)/n,$$
for $\alpha>3/4,~n=2$, or $\alpha>9/10,~n=3$.

The low frequency part in (5.1) is easy to control, and we can prove
$$\|T^*(\tilde{\psi}(\Delta)f)\|_{L^q(M)}\lesssim \|f\|_{L^2(M)}$$
for any $q\geq 2$. In fact, by making $t$ into a function $t(x)$, we only need to show
$$\|e^{-it(x)\Delta}\tilde{\psi}(\Delta)f\|_{L^q(M)}\lesssim \|f\|_{L^2(M)}.$$
By the compact support of $\tilde{\psi}$, we have
\begin{eqnarray*}
\|e^{-it(x)\Delta}\tilde{\psi}(\Delta)f\|_{L^q}&=& \left\|\sum_{\lambda_j\leq c_0}e^{-it(x)\lambda_j^2}\tilde{\psi}(\lambda_j^2)\hat{f}_j
e_j(x)\right\|_{L^q}\\
&\lesssim & \sum_{\lambda_j\leq c_0}\left\|e^{-it(x)\lambda_j^2}\tilde{\psi}(\lambda_j^2)\hat{f}_j
e_j(x)\right\|_{L^q}\\
&\lesssim & \sum_{\lambda_j\leq c_0}|\hat{f}_j|\left\|e_j(x)\right\|_{L^q}\\
&\lesssim & \sum_{\lambda_j\leq c_0}\lambda_j^{\delta(q)}|\hat{f}_j|\left\|e_j(x)\right\|_{L^2},
\end{eqnarray*}
for some positive $\delta (q)$. In the last step we applied the $L^q$ estimate for eigenfunctions of $-\Delta$ (see Theorem 5.1.1 of
\cite{SoFIO} or \cite{SoJFA}). If we take all the eigenfunctions to be $L^2$ normalized, the last term above is clearly
bounded by $\|f\|_{L^2}$ after using Schwartz's inequality.

To handle the rest of the terms in (5.1), we prove that for $0<h\leq1$,
$$\|T^*(\psi(h^2\Delta)f)\|_{L^p(M)}\lesssim h^{-\alpha}\|\psi(h^2\Delta)f\|_{L^2(M)},\eqno(5.2)$$
where $\alpha=3/4$ if $n=2$ or $\alpha=9/10$ if $n=3$. If (5.2) we prove, then
\begin{eqnarray*}
\sum_{k=1}^{+\infty}\|T^*(\psi(2^{-2k}\Delta)f)\|_{L^p}&\lesssim & \sum_{k=1}^{+\infty} 2^{\alpha k}\|\psi(2^{-2k}\Delta)f\|_{L^2}\\
&=&  \sum_{k=1}^{+\infty}2^{-\epsilon k} \|2^{(\alpha+\epsilon) k}\psi(2^{-2k}\Delta)f\|_{L^2}\\
&\lesssim & \left(\sum_{k=1}^{+\infty}2^{-2\epsilon k} \right)^{1/2}\left(\sum_{k=1}^{+\infty}\|(I-\Delta)^{(\alpha+\epsilon)/2}\psi(2^{-2k}\Delta)f\|^2_{L^2}\right)^{1/2}\\
&\leq & C_\epsilon \|(I-\Delta)^{(\alpha+\epsilon)/2}f\|_{L^2}\\
&=& C_\epsilon \|f\|_{H^{\alpha+\epsilon}}.
\end{eqnarray*}

Now we are left to prove (5.2). To do this, we need the following Strichartz estimate. \\\\
{\bf Lemma 5.1.} Let $0<h\leq 1$ and $p=\frac{2(n+2)}{n}$. Then
$$\|e^{-it\Delta}(\psi(h^2\Delta)f)\|_{L^p((0,1]\times M)}\lesssim h^{-1/p}\|\psi(h^2\Delta)f\|_{L^2(M)}.$$
\\
{\bf Proof.} This lemma can be inferred from the general Strichartz estimate in \cite{BGT}, Theorem 1, which was proved by
applying Keel-Tao's (\cite{KT}) theorem, after they constructed a parametrix for the frequency localized Schr\"{o}dinger equation
in local coordinate and proved a very short time version of the dispersion estimate. Here, we would like to present another proof which
applies H\"{o}rmander's oscillatory integral estimates. First we state the frequency localized parametrix as Lemma 2.7 of \cite{BGT}.\\\\
{\bf Parametrix.} Let $U_1$ be an open ball in $\Bbb R^n$ endowed with a Riemannian metric $g$. Take $U_2$ be an open ball in $U_1$, $\chi_0\in
C_0^\infty(U_2)$ and $\phi \in C_0^\infty(\Bbb R^n)$. Then for every $h\in(0,1]$ and $w_0\in C_0^\infty(U_1)$, there exists an $\alpha>0$
and $\widetilde{w}(s,x)\in C_0^\infty ([-\gamma,\gamma]\times U_2)$ which solves
$$\left\{\begin{array}{ll}
ih\partial_t\tilde{w}-h^2\Delta_g \tilde{w}=r,\\
\tilde{w}(0,x)=\chi_0(x)\phi(hD)w_0(x)\end{array}\right.$$
with $r(s,x)\in C_0^\infty ([-\gamma,\gamma]\times U_2)$ satisfying
$$\|r(s,x)\|_{L^\infty([-\gamma,\gamma],L^p(U_2))}\lesssim C_N h^N\|w_0\|_{L^2(U_1)},~~\forall~N.\eqno(5.3)$$
Furthermore, we have
$$\|\tilde{w}(s,x)\|_{L^p([-\gamma,\gamma]\times U_2)}\lesssim h^{\frac{n+1}{p}-\frac{n}{2}}\|w_0\|_{L^2(U_1)}.\eqno(5.4)$$
\\
Here we sketch the construction of $\widetilde{w}$ from \cite{BGT}. Consider
$$\widetilde{w}(s,x)=(2\pi h)^{-n}\int_{\Bbb R^n}e^{\frac{i}{h}\Phi(s,x,\xi)}a(s,x,\xi,h)\widehat{w}_0(\xi/h)d\xi\eqno(5.5)$$
where
$$a(s,x,\xi,h)=\sum_{j=0}^N h^ja_j(s,x,\xi)$$
Here $N$ is to be chosen large enough, $a_j\in C_0^\infty ([-\gamma,\gamma]\times U_2\times \Bbb R^n)$, with initial constrains
$$a_0(0,x,\xi)=\chi_0(x)\phi(\xi),~~a_j(0,x,\xi)=0,~~j\geq 1,$$
and $\Phi\in C^\infty ([-t_0,t_0]\times U_2\times B)$, where $B$ is a ball containing the support of $\phi$, with initial constrain
$\Phi(0,x,\xi)=x\cdot \xi$. Then the equations for $\phi$ and $a_j$ are given by the eikonal equation
$$\partial_s\Phi+\sum_{1\leq i,j\leq n}g^{ij}\partial_i\Phi\partial_j\Phi=0$$
and the transport equations
$$\partial_sa_0+2g(\nabla_g\Phi,\nabla_ga_0)+\Delta_g(\Phi)a_0=0$$
$$\partial_sa_j+2g(\nabla_g\Phi,\nabla_ga_j)+\Delta_g(\Phi)a_j=-\Delta_g(a_{j-1}),~~j\geq 1.$$
By the proof of Lemma 2.7 in \cite{BGT}, we also know that
$$r(s,x)=h^{N+2}(2\pi h)^{-n}\int_{\Bbb R^n}e^{\frac{i}{h}\Phi(s,x,\xi)}b(s,x,\xi,h)\widehat{w}_0(\xi/h)d\xi$$
for some $b\in C_0^\infty ([-\gamma,\gamma]\times U_2\times B)$. This easily yields (5.3).

Next we apply Lemma 2.1 to (5.5). The existence of phase function $\Phi$ on a small interval for $s$ is guaranteed by the Hamilton-Jacoby Theory. It is easy to see that the two conditions are satisfied when s=0. Actually, when $s=0$, for each $x$, $S(0,x)$ is a parabola. Then by continuity and compactness, the two conditions are satisfied for $s<\delta$ for some fixed $\delta=\delta(M)>0$.  Thus we get
$$\|\widetilde{w}(s,x)\|_{L^p([-\gamma,\gamma]\times U_2)}\lesssim h^{\frac{n+1}{p}-n}\|\widehat{w}_0(\xi/h)\|_{L^2(\Bbb R^n)} \lesssim h^{\frac{n+1}{p}-\frac{n}{2}}\|{w_0}(x)\|_{L^2(U_1)}.$$

Now let us continue to prove Lemma 5.1. Denote
$$\widetilde{w}_l(s,x)=e^{-ihs\Delta_g}(\chi_0\phi(hD)w_0)(s)$$
which is the solution to linear equation
$$\left\{\begin{array}{ll}
ih\partial_t\widetilde{w}-h^2\Delta \widetilde{w}=0,\\
\widetilde{w}(0,x)=\chi_0(x)\phi(hD)w_0(x).\end{array}\right.$$
Therefore,
$$\widetilde{w}(s,x)-\widetilde{w}_l(s,x)=\int_0^s e^{-ih(s-\tau)\Delta_g}r(\tau,x)d\tau.$$
So by (5.3), we have
$$\|\widetilde{w}(s,x)-\widetilde{w}_l(s,x)\|_{L^p([-\gamma,\gamma]\times U_2)}\lesssim h^N\|w_0\|_{L^2(U_1)}.$$
Thus (5.4) also holds for $\widetilde{w}_l(s,x)$. Note $\phi(hD)$ is a cut off on frequencies in $\Bbb R^n$. But in order to prove the lemma,
we need to show
$$\|e^{-ihs\Delta}(\psi(h^2\Delta)f)\|_{L^p([-\gamma,\gamma]\times M)}\lesssim h^{\frac{n+1}{p}-\frac{n}{2}}\|f\|_{L^2(M)}\eqno(5.6)$$
where the cut off $\psi(h^2\Delta)$ is made on frequencies $\lambda_j$. To treat this difference, we apply Corollary 2.4 of \cite{BGT}, which
actually says that there is a pseudodifferetial operator $\Psi(D)$ of order $0$ on $M$, such that, in local coordinates, $\Psi(\xi)$ is compactly
supported and
$$\|(I-\Psi(hD))\psi(h^2\Delta)f\|_{H^\sigma(M)}\lesssim C_{\sigma,N}h^N\|f\|_{L^2(M)}$$
holds for all $h\in (0,1]$, $\sigma>0,~N>0$ and $f\in C^\infty(M)$. Combining this with (5.4) for $\widetilde{w}_l$ and the boundedness of $e^{-ihs\Delta}$ on
$H^\sigma(M)$, we then reach (5.6) by constructing partitions of unity.

Finally, let us see how (5.6) implies Lemma 5.1. With a change of variable $hs\rightarrow t$, since $p=\frac{2(n+2)}{n}$, (5.6) implies
$$\|e^{-it\Delta}(\psi(h^2\Delta)f)\|_{L^p([-\gamma h,\gamma h]\times M)}\lesssim h^{\frac{n+2}{p}-\frac{n}{2}}\|f\|_{L^2(M)}=
\|f\|_{L^2(M)}$$
It is also easy to see that we can replace $f$ in the above $L^2$ norm by $\psi(h^2\Delta)f$ so that
$$\|e^{-it\Delta}(\psi(h^2\Delta)f)\|_{L^p([-\gamma h,\gamma h]\times M)}\lesssim \|\psi(h^2\Delta)f\|_{L^2(M)}.$$
Set $I_k=[(k-1)\gamma h,k\gamma h]$. Then
\begin{eqnarray*}
\|e^{-it\Delta}(\psi(h^2\Delta)f\|^p_{L^p((0,1]\times M)}& = & \sum_{k=1}^{(\gamma h)^{-1}}
\|e^{-it\Delta}\psi(h^2\Delta)f\|^p_{L^p(I_k\times M)}\\
&\lesssim & \sum_{k=1}^{(\gamma h)^{-1}}\|e^{-ik\gamma h\Delta}\psi(h^2\Delta)f\|^p_{L^2(M)}\\
&\lesssim & \sum_{k=1}^{(\gamma h)^{-1}}\|\psi(h^2\Delta)f\|^p_{L^2(M)}\\
&\lesssim & h^{-1}\|\psi(h^2\Delta)f\|^p_{L^2(M)},
\end{eqnarray*}
which proves Lemma 5.1. $\square$\\\\
{\bf Lemma 5.2.} Let $q\geq 2$. Suppose we have the following Strichartz estimate,
$$\|e^{-it\Delta}(\psi(h^2\Delta)f)\|_{L^q((0,1]\times M)}\lesssim h^{-\beta}\|\psi(h^2\Delta)f\|_{L^2(M)}.$$
Then for the maximal Schr\"{o}dinger operator $T^*$,  the following estimate holds,
$$\|T^*(\psi(h^2\Delta)f\|_{L^q(M)}\lesssim h^{-2/q-\beta}\|\psi(h^2\Delta)f\|_{L^2(M)}+\|\psi(h^2\Delta)f\|_{L^q(M)}.\eqno(5.7)$$
\\
{\bf Proof.} We will need the following inequality from \cite{Lee}, which says
$$\sup_{t\in [a,b]}|g(t)|\leq C_p\left(|g(a)|+\mu^{1/q-1}\|g'(t)\|_{L^q[a,b]}+\mu^{1/q}\|f\|_{L^q[a,b]}\right)\eqno(5.8)$$
for any smooth $g(t)$ on $[a,b]$, $\mu>0$ and $q\geq 1$.

Take $[a,b]=[0,1]$ and $g(t)=e^{-it\delta}\psi(h^2\Delta)f$. By (5.8) and Lemma 5.1,
\begin{eqnarray*}
\|T^*(\psi(h^2\Delta)f\|_{L^q(M)}&\lesssim & \|\psi(h^2\Delta)f\|_{L^q(M)}+\mu^{1/q} \|e^{-it\Delta}\psi(h^2\Delta)f\|_{L^q((0,1]\times M)}\\
&&+\mu^{1/q-1} \|e^{-it\Delta}(-i\Delta)\psi(h^2\Delta)f\|_{L^q((0,1]\times M)}\\
&\lesssim & \|\psi(h^2\Delta)f\|_{L^q(M)}+\mu^{1/q}h^{-\beta} \|\psi(h^2\Delta)f\|_{L^2}\\
&&+\mu^{1/q-1}h^{-\beta} \|(-i\Delta)\psi(h^2\Delta)f\|_{L^2}\\
&\lesssim & \|\psi(h^2\Delta)f\|_{L^q}+\mu^{1/q}(h^{-\beta}+\mu^{-1}h^{-2-\beta})\|\psi(h^2\Delta)f\|_{L^2}.
\end{eqnarray*}
By taking $\mu=h^{-2}$, we finish the proof of Lemma 5.2. $\square$\\

Now we are in a position to finish the proof of inequality (5.2), hence the main theorem. Let us combine the two lemmas together to get
$$\|T^*(\psi(h^2\Delta)f)\|_{L^p(M)}\lesssim h^{-3/p}\|\psi(h^2\Delta)f\|_{L^2(M)}+\|\psi(h^2\Delta)f\|_{L^p(M)}.$$
By Sobolev imbedding, the last term above is no larger than
$$\|\psi(h^2\Delta)f\|_{H^{n/2-n/p}(M)}\simeq h^{-(n/2-n/p)}\|\psi(h^2\Delta)f\|_{L^2(M)}.$$
Note $p=\frac{2(n+2)}{n}$. So
$$\frac{n}{2}-\frac{n}{p}=\frac{n}{n+2},~\frac{3}{p}=\frac{3n}{2(n+2)}.$$
A simple calculation then yields
$$\|T^*(\psi(h^2\Delta)f)\|_{L^p(M)}\lesssim h^{-3/4}\|\psi(h^2\Delta)f\|_{L^2(M)}$$
when $n=2$ and
$$\|T^*(\psi(h^2\Delta)f)\|_{L^p(M)}\lesssim h^{-9/10}\|\psi(h^2\Delta)f\|_{L^2(M)}$$
when $n=3$.
The case of $n=1$ will be dealt with in the next section.

\section{Solution on flat Torus and other special manifolds}
We may be able to improve our result if we could get better Strichartz estimate than in Lemma 5.1.
The argument in Section 5 comes from a trial to improve the Strichartz estimates \cite{BGT} on general manifolds. Although for general manifolds, we still get the same index and same loss, in some special manifolds, we do have more
precise Strichartz type inequalities, which enable us to get improved theorems.\\

To continue with the 2 dimensional flat Torus case, we will need the following Strichartz estimate on $\Bbb T^n$. It can be inferred from Bourgain's work \cite{B4}, Proposition 3.6.\\\\
{\bf Lemma 6.1.} For $n\geq 2$, the following Strichartz estimate holds
$$\|e^{-it\Delta}f\|_{L^4((0,1]\times \Bbb T^n)}\lesssim \|f\|_{H^s(\Bbb T^n)},~s>\frac{n}{4}-\frac{1}{2}.$$\\

Thus, by Lemma 5.2, we have
$$\|T^*(\psi(h^2\Delta)f)\|_{L^4(\Bbb T^2)}\lesssim h^{-1/2-s}\|\psi(h^2\Delta)f\|_{L^2(\Bbb T^2)}+\|\psi(h^2\Delta)f\|_{L^4(\Bbb T^2)}.$$
By Sobolev embedding
$$\|\psi(h^2\Delta)f\|_{L^4(\Bbb T^2)}\lesssim h^{-1/2}\|\psi(h^2\Delta)f\|_{L^2(\Bbb T^2)}.$$
Then
$$\|T^*(\psi(h^2\Delta)f)\|_{L^q(\Bbb T^2)}\lesssim h^{-1/2-s}\|\psi(h^2\Delta)f\|_{L^2(\Bbb T^2)}$$
for any $s>0$. By a similar argument as in Section 5, we finish the proof of the 2-dimensional flat Torus case.

For the $n=1$ case of Theorem 1.2, we first notice that all connected, compact 1 dimensional manifolds are isometric to circles. So we only need to consider $\Bbb T^1$. We have the following Strichartz estimate from \cite{B4}, Proposition 2.36.\\\\
{\bf Lemma 6.2.} The following Strichartz estimate holds
 $$\|e^{-it\Delta}f\|_{L^6((0,1]\times \Bbb T^1)}\lesssim \|f\|_{H^s(\Bbb T^1)},~s>0.$$\\
As above, we conclude
$$\|T^*(\psi(h^2\Delta)f)\|_{L^6(\Bbb T^1)}\lesssim h^{-1/3-s}\|\psi(h^2\Delta)f\|_{L^2(\Bbb T^1)}.$$
Thus we finish the proof of Theorem 1.2 for $n=1$.

Now let us consider the higher dimensional flat torus $\Bbb T^n$, $n\geq 3$. By applying the stronger Strichartz estimate in the following Lemma and the argument above, we can reduce the amount of regularity requirement to some number less than 1 for all dimensional flat torus. \\\\
{\bf Lemma 6.3.} The following Strichartz estimate holds
$$\|e^{-it\Delta}f\|_{L^q((0,1]\times \Bbb T^n)}\lesssim \|f\|_{H^s(\Bbb T^n)},~s>0,~q\leq\frac{2(n+1)}{n}.$$\\
The proof of this lemma can be found in \cite{JB}. As a consequence, we have the following improved theorem.\\\\
{\bf Theorem 6.1.} Let $e^{-it\Delta}$ be the Schr\"{o}dinger operator defined on $\Bbb T^n$. Then $e^{-it\Delta}f$ converges pointwisely to $f$ if $f\in H^\alpha(\Bbb T^n)$, where $\alpha>\frac{n}{n+1}$, $n\geq 3$.\\

Finally let us consider one type of manifold whose geodesics are closed with a common period. For the geometric properties of such manifolds, one is referred to \cite{Be}.
Here we only apply the Strichartz estimate (\cite{BGT}, Theorem 4) for the Schr\"{o}dinger operator on such manifolds,
$$\|e^{-it\Delta}f\|_{L^4((0,1]\times M)}\lesssim \|f\|_{H^s(M)},~s>s_0(n)$$
where $s_0(2)=1/8$, $s_0(n)=n/4-1/2$ for $n\geq 3$. Note the $n$ sphere $\Bbb S^n$ is one of the above manifolds. Furthermore, the loss $s_0(n)$ has been proved to be sharp for $\Bbb S^n$. Similarly, we have\\\\
{\bf Theorem 6.2.} Let $M$ be the manifold described above. Then $e^{-it\Delta}f$ converges pointwisely to $f$ if\\
(i)~~$n=2$ and $f\in H^\alpha(\Bbb M^2),~\alpha>5/8$, or\\
(ii) $n=3$ and $f\in H^\alpha(\Bbb M^3),~\alpha>3/4$.\\

{\bf Acknowledgement} The authors are grateful for helpful suggestions and comments from C.D. Sogge, S. Lee and C.B. Wang.


\begin{thebibliography}{99}

\bibitem{Be} A. Besse, Manifolds all of whose geodesics are closed, Springer-Verlag, New York, 1978.

\bibitem{B1} J. Bourgain, A remark on Schr\"{o}dinger operators, Israel J. Math. 77 (1992), 1-16.

\bibitem{B2} J. Bourgain, On the Schr\"{o}dinger maximal founction in higher dimension,  J. Proc. Steklov Inst. Math. 280 (2013), 46-60

\bibitem{B3} J. Bourgain, Exponential sums and nonlinear Schr\"{o}dinger equations, Geom. Funct. Anal. 3 (1993), 157-178.

\bibitem{B4} J. Bourgain, Fourier transform restriction phenomena for certain lattice subsets and applications
to nonlinear evolution equations, I. Schr\"{o}dinger equations, Geom. Funct. Anal. 3 (1993), 107-156.

\bibitem{JB}J. Bourgain, Moment inequalities for trigonometric polynomials with spectrum in curved hypersurfaces, Israel Journal of Mathematics, 193 (2013), 441-458.

\bibitem{BGT} N. Burq, P. G\'{e}rard and N. Tzvetkov, Strichartz inequalities and the nonlinear Schr\"{o}dinger
equation on compact manifolds, Amer. J. Math. 126 (2004), 569-605.

\bibitem{Car}A. Carbery, Radial Fourier multipliers and associated
maximal functions, Recent progress in Fourier analysis,
North-HollandMath. Stud. 111, North-Holland, Amsterdam, 1985,
5-45.

\bibitem{Ca}L. Carleson, Some Analytic problems related to Statistical Mechanics, Euclidean harmonic
analysis, Lecture Notes in Math. 779, Springer, Berlin, 1980,
5-45.

\bibitem{Cow}M. Cowling, Pointwise behavior of solutions to Schr\"{o}dinger equations, Harmonic analysis,
Lecture Notes in Math. 992, Springer, Berlin, 1983, 83-90.

\bibitem{DK}B.E.J. Dahlberg and C.E. Kenig, A note on the almost everywhere hehavlor of solutions to the
Schr\"{o}dinger equation, Lecture Notes in Math. 908,
Springer-Verlag, Berlin and New York, 1982, 205-208.


\bibitem{Ho}L. H\"{o}mander, The Analysis of Linear Partial Differential Operators III: Pseudo-Differential Operators, Springer-Verlag Berlin Heidelberg 2007.

\bibitem{KR}C.E. Kenig and A. Ruiz, A strong type (2,2) estimate for a maximal operator associated to the
Schr\"{o}dinger equation, Trans. Amer. Math. Soc. 280 (1983), 239-246.

\bibitem{KT}M. Keel and T. Tao, End point Strichartz estimates, Amer. J. of Math. 120 (1998), 955-980.

\bibitem{Lee}S. Lee, On pointwise convergence of the solutions to Schr\"{o}dinger equations in $\Bbb R^2$, Int. Math.
Res. Not. 2006, Art. ID 32597, 1-21.

\bibitem{LR} R. Luc?a and K. Rogers, An improved necessary condition for the Schr\"{o}dinger maximal estimate, ArXiv: 1506.05325v1.

\bibitem{MVV} A. Moyua, A. Vargas and L. Vega, Schr\"{o}dinger maximal function and restriction properties
of the Fourier transform, Int. Math. Res. Not. 16 (1996), 793-815.

\bibitem{MV} Moyua, A.(E-BILB); Vega, L.(E-BILB), Bounds for the maximal function associated to periodic solutions of one-dimensional dispersive equations. Bull. Lond. Math. Soc. 40 (2008), no. 1, 117?128. 

\bibitem{SD1} S. Doi, Smoothing effects of Schr\"{o}dinger evolution groups on Riemannian manifolds.
Duke Math. J. 82 (1996), 679-706.

\bibitem{SD2} S. Doi, Smoothing effects for Schr\"{o}dinger evolution equation and global behavior of geodesic flow.
Math. Ann. 318 (2000), 355-389.

\bibitem{Sj}P. Sj\"{o}lin, Regularity of solutions to the Schr\"{o}dinger equation, Duke Math. J. 55 (1987),
699-715.

\bibitem{SoHZ}C.D. Sogge, Hangzhou lectures on eigenfunctions of the Laplacian, Annals of Mathematics Studies, 188 Princeton University Press, Princeton, NJ, 2014.


\bibitem{SoFIO}C.D. Sogge, Fourier integrals in classical analysis. Cambridge Tracts in Mathematics, 105. Cambridge University Press, Cambridge, 1993.

\bibitem{SoJFA}C.D. Sogge, Concerning the $L^p$ norm of spectral clusters for second-order elliptic operators on compact manifolds,
J. Funct. Anal. 77 (1988), 123-138.

\bibitem{T}T. Tao, A sharp bilinear restrictions estimate for paraboloids, Geom. Funct. Anal. 13 (2003),
1359-1384.

\bibitem{TV1}T. Tao and A. Vargas, A bilinear approach to cone multipliers. I. Restriction estimates,
Geom. Funct. Anal. 10 (2000), 185-215.

\bibitem{TV2}T. Tao, A. Vargas, A bilinear approach to cone multipliers. II. Applications, Geometric and Functional Analysis 10 (2000), 216-258.

\bibitem{Ve1}L. Vega, Schr\"{o}dinger equations: pointwise convergence to the initial data, Proc. Amer. Math.
Soc. 102 (1988), 874-878.

\bibitem{BG} B.G. Walther, Some $L^p(L^1)$- and $L^2(L^2)$- estimates for oscillatory Fourier transforms,
Analysis of Divergence (Orono, ME, 1997), 213-231, Appl. Numer. Harmon. Anal., Birkh¡§aser Boston, Boston, MA, 1999.





\end{thebibliography}
\end{document}